\documentclass[12pt]{article}
\usepackage{hyperref}
\usepackage{amsmath}
\usepackage{mathtools}
\usepackage{epsf}
\usepackage{theorem}
\usepackage{indentfirst}
\usepackage{latexsym}
\usepackage{amssymb}
\usepackage[dvips]{graphicx}
\usepackage{flafter}
\usepackage{caption}
\usepackage{textcomp}
\usepackage{supertabular}
\usepackage{longtable, booktabs}
\usepackage{hhline}
\usepackage{bigstrut,bigdelim}
\usepackage{rotating}
\usepackage{multicol}
\usepackage{multirow}
\usepackage{makeidx}
\usepackage{epsfig}
\usepackage{fancyhdr}
\newif\ifger
\gerfalse
\newtheorem{theorem}{Theorem}
\newtheorem{lemma}{Lemma}[section]

\newtheorem{remark}{Remark}
\begin{document}
\baselineskip=19pt
\title{Block-transitive and point-primitive $2$-$(v,k,2)$ designs with
 sporadic socle
 }
\author{Xiaohong Zhang, Shenglin Zhou\footnote{Corresponding author. This work is
supported by the National Natural Science Foundation
of China (Grant No.11471123). slzhou@scut.edu.cn}\\
\small \it School of Mathematics, South China University of Technology,\\
\small \it Guangzhou 510641, P.R. China}

\date{}
\maketitle

\begin{abstract}
   The purpose of this paper is to classify all pairs $(\mathcal{D}, G)$, where $\mathcal{D}$ is a non-trivial $2$-$(v, k, 2)$ design, and $G\leq Aut(\mathcal{D})$ acts transitively on the set of blocks of $\mathcal{D}$ and primitively on the set of points of $\mathcal{D}$ with sporadic socle. We prove that there exists only one such pair $(\mathcal{D}, G)$ in which $\mathcal{D}$ is a $2$-$(176,8,2)$ design and $G=HS$, the Higman-Sims simple group.

\medskip
\noindent{\bf MR(2000) Subject Classification:} 05B05, 05B25,
20B25

\end{abstract}

\section{Introduction}

  A $2$-$(v, k, \lambda)$ {\it  design} $\mathcal{D}$ is a set of $v$ points $\mathcal{P}$ together with a collection $\mathcal{B}$ of distinct $k$-subsets of $\mathcal{P}$, called blocks, such that any two points lie in exactly $\lambda$ blocks. We denote the number of blocks by $b$ and the number of blocks containing a point by $r$. We shall also assume that $\mathcal{P}$ is finite, $b>1$ and $k>2$. It is well known that
\begin{align*}
 bk&=vr;\\
\lambda(v-1)&=r(k-1);\\
b &\geq v.
 \end{align*}

An {\it automorphism} of $\mathcal{D}$ is a permutation of the points of $\mathcal{D}$ that also permutes the blocks of $\mathcal{D}$. The {\it full automorphism group}  of $\mathcal{D}$ is the group of all automorphisms of $\mathcal{D}$ and is denoted by $Aut(\mathcal{D})$.
If $G \leq Aut(\mathcal{D})$, then $G$ is called an {\it automorphism group} of $\mathcal{D}$. We say that $G$ is {\it block-transitive} (resp. {\it point-transitive}) if $G$ acts transitively on $\mathcal{B}$ (resp. $\mathcal{P}$) and {\it point-primitive} if $G$ acts primitively on $\mathcal{P}$.

By a {\it flag} of $\mathcal{D}$ we mean an incident pair $(\alpha, B)$, where $\alpha$ is a point of $\mathcal{D}$ and $B$
is a block of $\mathcal{D}$. If $G \leq Aut(\mathcal{D})$, then $G$ is {\it flag-transitive} if it is transitive on the set of flags of $\mathcal{D}$. An {\it antiflag} of $\mathcal{D}$ is a non-incident pair $(\alpha, B)$, where $\alpha$ is a point of $\mathcal{D}$ and $B$
is a block of $\mathcal{D}$. If $G \leq Aut(\mathcal{D})$, then $G$ is {\it antiflag-transitive} if it is transitive on the set of antiflags of $\mathcal{D}$.

The classification of block-transitive $2$-$(v,k,1)$ designs is now under way and plentiful results have been achieved, see \cite{CaGiZa,CaNePr,CaPr}. Here, we just mention some of results which have been obtained for automorphism with small $k$. In \cite{Camina}, the classification of $2$-$(v,4,1)$ designs with a block-transitive and solvable group of automorphisms has been achieved by Camina and Siemons, and in \cite{Li} Li completed the follow-up unsolvable case. In \cite{AR}, Camina and Mischke classified block-transitive and point-imprimitive $2$-$(v,k,1)$ designs with $k\leq 8$. In 2009, Betten, Delandtsheer, Law et al. \cite{B} further classified this type of $2$-$(v,k,1)$ designs satisfying $(k,v-1) \leq 8$. In particular, Camina and Spiezia proved in \cite{F} that if $G$ is an almost simple group which acts block-transitively on a $2$-$(v,k,1)$ design then $Soc(G)$ cannot be a sporadic group.

In the case of $\lambda=2$, the classification of flag-transitive symmetric $2$-$(v,k,2)$ designs has almost been completed by Regueiro in a sequence of 4 papers, see \cite{Re05a,Re05b,Re08a,Re08b}. Recently,
Liang and Zhou \cite{Liang} proved that a non-trivial non-symmetric $2$-$(v,k,2)$ design admitting a flag-transitive and point-primitive almost simple automorphism group $G$ with sporadic socle must be the unique $2$-$(176,8,2)$ design with $G=HS$, the Higman-Sims simple group. However, for block-transitive $2$-$(v,k,2)$ designs, there are only a few known results. Classifying block-transitive $2$-$(v,k,2)$ designs seems to be a challenging problem. In this paper, we give a complete classification of block-transitive and point-primitive $2$-$(v,k,2)$ designs with sporadic socle.

Our main result is the following:

\begin{theorem}\label{thm1}\quad
 Let $\mathcal{D}$ be a non-trivial $2$-$(v,k,2)$ design. Assume that $\mathcal{D}$ has an automorphism group $G$ that is block-transitive and point-primitive with $ Soc(G)$ a sporadic simple group. Then $\mathcal{D}$ must be the unique $2$-$(176,8,2)$ design with $G=HS$, the Higman-Sims simple group.

 \end{theorem}

\begin{remark}
{\rm

 \ \ \ As a matter of fact, the automorphism group $G=HS$ acting on the $2$-$(176,8,2)$ design is not only
   block-transitive and point-primitive, but also flag-transitive $\cite{Liang}$ and antiflag-transitive.}
\end{remark}

 \section{ Preliminaries }
In this section we collect some basic results that will be used throughout the proof of the main theorem.

\begin{lemma}\label{LEM1} $\cite{Block}$ \quad
 Let $\mathcal{D}$ be a non-trivial $2$-$(v,k,\lambda)$ design. If $G \leq$ $Aut(\mathcal{D})$ acts block-transitively on $\mathcal{D}$, then $G$ acts point-transitively on $\mathcal{D}$. \end{lemma}

\begin{lemma} \label{LEM2} \quad
 Let $\mathcal{D}$ be a $2$-$(v,k,2)$ design with $b$ blocks. Given an integer $v > 2$, there are only finitely many pairs $(k$, $b)$ of integers such that there exists a $2$-$(v,k,2)$ design with $b$ blocks admitting a block-transitive and point-primitive automorphism group $G$ with the following properties:
 \begin{enumerate}
 \item[$(1)$] \,
  $k-1\mid 2(v-1)$;
 \item[$(2)$] \,
  $2<k<v$;
 \item[$(3)$] \,
  $b=\frac{2v(v-1)}{k(k-1)}\in \mathbb{N}$;
 \item[$(4)$] \,
  $v\leq b$;
 \item[$(5)$] \,
  $b\mid  |$G$|$.
 \end{enumerate}
 \end{lemma}
\textbf{Proof.}
Parts \rm (1)-\rm (4) are general properties of designs. Part \rm (5) follows from $b=|G:G_B|$ for $B\in \mathcal{B}$ as $G$ is transitive on the set of blocks of $\mathcal{D}$.$\hfill\square$

The following lemma is crucial for the proof of Theorem \ref{thm1}.
\begin{lemma}  \label{LEM3}  \quad
Let $\mathcal{D}$ be a non-trivial $2$-$(v,k,2)$ design, and let $G \leq$ $Aut(\mathcal{D})$. Let $B$ be a block of $\mathcal{D}$. Then the block-length $k$ can be written as the sum of some orbit-lengths of $G_B$ on $\mathcal{P}$.

\end{lemma}
\textbf{Proof.}
Let $\Delta_1, \Delta_2,\ldots, \Delta_s$ be the orbits of $G_B$ on $\mathcal{P}$, so that $\mathcal{P}=\Delta_{1}\cup\Delta_2\cup\cdots\cup\Delta_s$. Since $G_B$ also acts on $B$, $B$ is an union of $G_B$-orbits, that is, $B=\Delta_{i_{1}}\cup\Delta_{i_{2}}\cup\cdots\cup\Delta_{i_{t}}$ for some $i_1, i_2,\ldots, i_t$ $\in$ $\{1, 2,\ldots, s\}$.
Consequently, the block-length $k$ is the sum of some orbit-lengths of $G_B$ on $\mathcal{P}$.$\hfill\square$

 \section{Proof of Theorem 1}

We prove Theorem \ref{thm1} in two subsections. In the first subsection, we apply the properties \rm (1)-\rm (5) in Lemma \ref{LEM2} and obtain a number of possible parameters $(v,b,r,k,\lambda)$ where $\lambda=2$. In the second subsection, we further analyse these potential designs to exclude or construct them.

 \subsection{Potential $2$-$(v,k,2)$ designs}\label{Sub3.1}

First we describe briefly how to search for potential block-transitive and point-primitive $2$-$(v,k,2)$ designs with sporadic socle.

 In this paper we restrict ourselves to the case when $Soc(G)$ is a sporadic simple group. Therefore we know that $G$ must be almost simple, that is $Soc(G)\unlhd G\leq Aut(Soc(G))$. Let $S$ be an arbitrary sporadic simple group. Then $|Out(S)| = 1$ or 2 from the {\sc Atlas} of Finite Group Representations $\cite{Con}$, which we will always refer to as `the {\sc Atlas}'. Since $S$ is a non-abelian simple group, it is clear that $Z(S)$, the center of $S$, has to be the identity. Combining the fact $S/Z(S) \cong Inn(S)$ with $Aut(S)/Inn(S)\cong Out(S)$, we obtain $S\unlhd G\leq S: Out(S)$, so $G=S$ or $S:2$.

Since $G$ is point-transitive (Lemma \ref{LEM1}), $G$ contains a subgroup $G_{\alpha}$, the stabilizer of a point $\alpha$, with index $v$. Since $G$ is block-transitive, $G$ contains a subgroup $G_B$, the stabilizer of a block $B$, with index $b$. As $G$ is point-primitive, the subgroup $G_{\alpha}$ must be maximal in $G$. For all sporadic groups, except the Monster, the complete list of maximal subgroups can be learned by consulting the {\sc Atlas}. Therefore, for each sporadic group, we can find the possible values for $|G_{\alpha}|$, and consequently, for $v$. For a fixed $v$, if no $k$ and $b$ satisfying properties \rm (1)-\rm (5) in Lemma \ref{LEM2} can be found, then we can exclude this value of $v$. Otherwise, for potential values of $k$ and $b$, we will construct a $2$-$(v,k,2)$ design with $b$ blocks or prove such a design does not exist.

  In the first part of the proof, apart from the Monster, we have examined all other 25 sporadic simple groups with the aid of the computer algebra system {\sf GAP} ($\cite{Gap}$). All possible parameters $(v,b,r,k,2)$ satisfying properties \rm (1)-\rm (5) are listed in Table 1 below. The case of the Monster will be dealt with separately in the next subsection.
 \bigskip
 \bigskip
 \begin{longtable}{clllc}
\caption{Potential $2$-designs and automorphism groups}\\ \hline
\endfirsthead
\multicolumn{5}{l}{(Continued)}\\
\hline
\endhead
\hline
\multicolumn{5}{r}{(Continued on next page)}\\
\endfoot\hline
\endlastfoot
{Case} & $G$ & $ G_{\alpha}$ &  $(v,b,r,k,\lambda)$ & {\small Lemma}\\
\hline
  $1$ & $M_{11}$ & $M_{10}$ & $(11,11,5,5,2)$  &\ref{LEM3.2} \\
  $2$ & $$ & $L_{2}(11)$ & $(12,44,11,3,2)$  &\ref{LEM3.1} \\
  $3$ & $$ & $M_{9}:2$ & $(55,990,54,3,2)$  &\ref{LEM3.3} \\
  $4$ &  &  & $(55,495,36,4,2)$ &\ref{LEM3.3}\\
  $5$ &  &  & $(55,66,12,10,2)$ &\ref{LEM3.4}\\
  $6$ & $M_{12}$ & $M_{11}$  &  $(12,44,11,3,2)$ &\ref{LEM3.1}\\
  $7$ & $M_{22}$ & $L_{3}(4)$& $(22,154,21,3,2)$ &\ref{LEM3.1}\\
  $8$ &  &  & $(22,77,14,4,2)$ &\ref{LEM3.2}\\
  $9$ &  &  & $(22,22,7,7,2)$ &\ref{LEM3.2}\\
  $10$ & $$ & $A_{7}$& $(176,560,35,11,2)$ &\ref{LEM3.1}\\
  $11$ & $M_{22}:2$ & $L_{3}(4):2_{2}$  & $(22,154,21,3,2)$ &\ref{LEM3.2}\\
  $12$ &  &  & $(22,77,14,4,2)$&\ref{LEM3.2}\\
  $13$ &  &  & $(22,22,7,7,2)$&\ref{LEM3.2}\\
  $14$ & $M_{23}$ & $L_{3}(4):2_{2}$ & $(253,21252,252,3,2)$ &\ref{LEM3.3}\\
  $15$ & &  &  $(253,10626,168,4,2)$&\ref{LEM3.2}\\
  $16$ & &  &  $(253,3036,84,7,2)$&\ref{LEM3.1}\\
  $17$ & &  &  $(253,2277,72,8,2)$&\ref{LEM3.1}\\
  $18$ & &  &  $(253,1771,63,9,2)$&\ref{LEM3.2}\\
  $19$ & &  &  $(253,276,24,22,2)$&\ref{LEM3.1}\\
  $20$ & $$ & $M_{11}$& $(1288,18216,198,14,2)$ &\ref{LEM3.1}\\
  $21$ & $$ & $2^4:(3\times A_{5}):2$& $(1771,17711,60,60,2)$ &\ref{LEM3.4}\\
  $22$ & $M_{24}$ & $M_{23}$ & $(24,184,23,3,2)$ &\ref{LEM3.1}\\
  $23$ & $$ & $M_{22}:2$ & $(276,5060,110,6,2)$ &\ref{LEM3.1}\\
  $24$ &  &  &  $(276,1380,55,11,2)$ &\ref{LEM3.1}\\
  $25$ & $$ & $M_{12}:2$ & $(1288,18216,198,14,2)$ &\ref{LEM3.1}\\
  $26$ & $$ & $2^6:3.S_{6}$ & $(1771,1771,60,60,2)$ &\ref{LEM3.2}\\
  $27$ & $HS$ & $M_{22}$ & $(100,3300,99,3,2)$ &\ref{LEM3.1}\\
  $28$ & &  &  $(100,1650,66,4,2)$&\ref{LEM3.1}\\
  $29$ & &  &  $(100,220,22,10,2)$&\ref{LEM3.1}\\
  $30$ & &  &  $(100,150,18,12,2)$&\ref{LEM3.1}\\
  $31$ & $$ & $U_{3}(5):2$ & $(176,1100,50,8,2)$ &$\mathcal{D}$\\
  $32$ &  & &  $(176,560,35,11,2)$ &\ref{LEM3.1}\\
  $33$ & $$ & $2\times A_{6}.2^2$ & $(15400,138600,531,59,2)$ &\ref{LEM3.2}\\
  $34$ & $HS:2$ & $M_{22}:2$ & $(100,3300,99,3,2)$ &\ref{LEM3.1}\\
  $35$ & &  &  $(100,1650,66,4,2)$&\ref{LEM3.1}\\
  $36$ & &  &  $(100,220,22,10,2)$&\ref{LEM3.1}\\
  $37$ & &  &  $(100,150,18,12,2)$&\ref{LEM3.1}\\
  $38$ & $$ & $(2\times A_{6}.2^2).2$ & $(15400,138600,531,59,2)$ &\ref{LEM3.2}\\
  $39$ & $J_{2}$ & $U_{3}(3)$  & $(100,150,18,12,2)$ &\ref{LEM3.1}\\
  $40$ & $J_{2}:2$ & $U_{3}(3):2$  & $(100,150,18,12,2)$ &\ref{LEM3.1}\\
  $41$ & $Co_{3}$ & $McL:2$ &  $(276,25300,275,3,2)$ &\ref{LEM3.1}\\
  $42$ & &  &  $(276,5060,110,6,2)$ &\ref{LEM3.1}\\
  $43$ & &  &  $(276,1380,55,11,2)$ &\ref{LEM3.1}\\
  $44$ & &  &  $(276,1150,50,12,2)$ &\ref{LEM3.1}\\
  $45$ & &  &  $(276,300,25,23,2)$ &\ref{LEM3.1}\\
  $46$ & $McL$ & $M_{22}$ &  $(2025,16200,184,23,2)$ &\ref{LEM3.1}\\
  $47$ & &  &  $(2025,14850,176,24,2)$ &\ref{LEM3.1}\\
  $48$ & $$ & $3^{1+4}:2.S_{5}$ &  $(15400,138600,531,59,2)$ &\ref{LEM3.1}\\
  $49$ & $$ & $5^{1+2}:3:8$ &  $(299376,779625,1250,480,2)$ &\ref{LEM3.2}\\
  $50$ & $McL:2$ & $3^{1+4}:4.S_{5}$ &  $(15400,138600,531,59,2)$ &\ref{LEM3.1}\\
  $51$ & $$ & $5^{1+2}:3:8.2$ &  $(299376,779625,1250,480,2)$ &\ref{LEM3.6}\\
  $52$ & $Fi_{24}'$ & $Fi_{23}$ &  $(306936,1904952,1955,315,2)$ &\ref{LEM3.5}\\
  $53$ & $Fi_{24}$ & $Fi_{23}\times 2$ &  $(306936,1904952,1955,315,2)$ &\ref{LEM3.5}\\
  $54$ & $J_{1}$ & $F_{168}$  & $(1045,2508,72,30,2)$ &\ref{LEM3.1}\\
  $55$ & $$ & $F_{114}$ &  $(1540,6270,114,28,2)$ &\ref{LEM3.1}\\
  $56$ & &  &  $(1540,1596,57,55,2)$ &\ref{LEM3.4}\\
  $57$ & $$ & $F_{110}$ &  $(1596,5852,110,30,2)$ &\ref{LEM3.4}\\
  $58$ & $$ & $D_{6}\times D_{10}$ &  $(2926,3990,90,66,2)$ &\ref{LEM3.5}\\
  $59$ & $O'N$ & $L_{3}(7):2$  & $(122760,1484280,1729,143,2)$ &\ref{LEM3.5}\\
  $60$ & $J_{3}$ & $3^{2+1+2}:8$  & $(25840,174420,594,88,2)$ &\ref{LEM3.4}\\
  $61$ & $J_{3}:2$ & $3^{2+1+2}:8.2$  & $(25840,174420,594,88,2)$ &\ref{LEM3.6}\\

\end{longtable}
\begin{remark}
{\rm

 \ \ \ In each case, the last column of Table 1 indicates that we rule out it by the lemma in subsection \ref{Sub3.2}, and the unique symbol $\mathcal{D}$ refers to the design we construct in Lemma $3. 8$ in the next subsection.}
\end{remark}

 \subsection{Analyzing parameters and corresponding groups}\label{Sub3.2}

In this section, we analyze the potential parameters in Table 1 and the corresponding automorphism groups. First, Lemmas 3.1-3.6 below prove that in Table 1 each case but the $\rm31^{st}$ cannot occur. Next, Lemma 3.7 deals with the Monster $M$, for which the list of maximal subgroups has not been finished yet. Finally, Lemma 3.8 deals with case 31 for which we obtain a unique design. This will complete the proof of Theorem \ref{thm1}.

The commands mentioned in the proof below are performed by the computer algebra system {\sc Magma} \cite{Bos}.

\begin{lemma}\label{LEM3.1}{\rm}\quad
Cases 2, 6, 7, 10, 16, 17, 19, 20, 22-25, 27-30, 32, 34-37, 39-48, 50, 54 and 55 cannot occur.\end{lemma}
\textbf{Proof.}
Since $G$ is block-transitive, the stabilizer $G_B$ of a block $B$ satisfies $|G:G_B|=b$. Applying the command ${\tt {Subgroups(G:OrderEqual:=n)}}$ where $n=|G|/b$ to each case, we find that such a subgroup $G_B$ does not exist.$\hfill\square$

\begin{lemma}\label{LEM3.2}{\rm}\quad
Cases 1, 8, 9, 11-13, 15, 18, 26, 33, 38 and 49 cannot occur.\end{lemma}
\textbf{Proof.}
 As an example, we analyze case $49$, as other cases are similar. Suppose $G=McL$. If $b=779625$, then $|G_B|=1152$. We then run the command ${\tt {Subgroups(G:OrderEqual:=n)}}$ where $n=1152$ and obtain two conjugacy  classes of subgroups with index 1152, denoted by $L_1$ and $L_2$ as representatives. If we run the commands ${\tt {O:=Orbits(L)}}$ for $L=L_1$ or $L_2$, and ${\tt {\sharp O[j]}}$, $j=1, 2, \ldots, 16$, we can see the 7 smallest orbit-lengths of $L_1$ are 144, $288^{5}$ and 576, and the 16 smallest orbit-lengths of $L_2$ are 48, $288^{4}$, $384^{10}$ and 576, where $a^{b}$ means that the orbit-length $a$ appears $b$ times. It is clear that there is no way to express $k=480$ as the sum of some orbit-lengths of $L_1$ or $L_2$. Hence a $2$-$(v,k,2)$ design with $b$ blocks admitting a block-transitive and point-primitive automorphism group $G$ dose not exist by Lemma \ref{LEM3}.$\hfill\square$

\begin{lemma}\label{LEM3.3}{\rm}\quad
Cases 3, 4 and 14 cannot occur.\end{lemma}
\textbf{Proof.}
 In case 4, $G=M_{11}$ has only one conjugacy class of subgroups of index $b=|G:G_B|$, where $B\in\mathcal{B}$, and $G_B$ has seven orbits $\Delta_1, \Delta_2,\ldots, \Delta_7$ on $\mathcal{P}$ with lengths 1, 2, 4, $8^{2}$, $16^{2}$, respectively. By Lemma 2.3, $B$ is union of $G_B$-orbits. For $k=4$, we have $B=\Delta_3$ and so $\Delta_3$ is a block. However, the block orbit-length $|\Delta_3^{G}|=165\neq495$, contradicting the fact that $G$ is block-transitive.

 Cases 3 and 14 can be ruled out similarly.$\hfill\square$

\begin{lemma}\label{LEM3.4}{\rm}\quad
Cases 5, 21, 56, 57 and 60 cannot occur.\end{lemma}
\textbf{Proof.}
In case 5, $G=M_{11}$ has only one conjugacy class of subgroups of index $b=|G:G_B|$, where $B\in\mathcal{B}$, and $G_B$ has three orbits $\Delta_1$, $\Delta_2$, $\Delta_3$ on $\mathcal{P}$ with lengths 10, 15, 30, respectively. Since $B$ is union of $G_B$-orbits and $k=10$, we have $B=\Delta_1$ and so $\Delta_1$ is a block. It is easy to work out the block orbit-length $|\Delta_1^{G}|=66=b$. However, the command ${\tt {Design<2,55|C>}}$ where $C=\Delta_1^{G}$ returns that the structure is not a $2$-design, a contradiction.
The other cases can be ruled out similarly.$\hfill\square$

\begin{lemma}\label{LEM3.5}{\rm}\quad
Cases 52, 53, 58 and 59 cannot occur.\end{lemma}
\textbf{Proof.}
In case 52, suppose $G=Fi_{24}'$ and $b=1904952$, so that $|G_B|=$ 658917237384806400. However, in $Fi_{24}'$, there is no maximal subgroup whose order is divisible by $|G_B|$ (see  \cite{Con}). This means that $G_B$ is contained in none of the maximal subgroups of $G$, contradicting the definition of a maximal subgroup (see $\cite{Dem}$).

Similarly, we can prove that cases 53 and 59 cannot occur.

In case 58, suppose $G=J_{1}$ and $b=3990$, so that $G_B$ has order 44. Inspecting the list of maximal subgroups of $J_{1}$ and their orders (see \cite{Con}), the only possibility is that $G_B$ is contained in a maximal subgroup of $J_{1}$ of order 660. This group is isomorphic to $L_{2}(11)$, but $L_{2}(11)$ does not contain any maximal subgroup of order divisible by 44 (see \cite{Con}), a contradiction.$\hfill\square$

\begin{lemma}\label{LEM3.6}{\rm}\quad
Cases 51 and 61 cannot occur.\end{lemma}
\textbf{Proof.}
In case 51, with the help of the command ${\tt {M:=PermRepKeys(McL:2)}}$, we can get the permutation representations of the group $McL:2$ of degrees 275, 4050, 7128, 22275, 44550, respectively. Since the degree needed is 299376, we will construct the permutation representation of the group $McL:2$ of degree 299376. Let $G$ be the fifth permutation representation of the group $McL:2$, which is of degree 44550. Of course, we can also choose permutation representation of other degrees, which will lead to the identical results. It is easy to find all maximal subgroups of $G$, and we only choose one, say $H$, with order $|G|/299376=6000$. Now by running the command ${\tt {F1,N:=CosetAction(G,H)}}$, we successfully construct a group $N$, which is the very permutation representation of the group $McL:2$ of degree 299376. We find that 299376 is too big to determine the subgroups of $N$ with index $b=779625$, but 44550 is appropriate to determine the subgroups of $G$ with the same index. So in consideration of the isomorphism between $G$ and $N$, we decide to construct a map between generators of these two groups. Since both $G$ and $N$ have two generators, the command ${\tt {F:=hom< G -> N | G.1-> N.1,G.2-> N.2 >}}$ gives a map from generators of $G$ to generators of $N$. We can easily get the subgroups of $G$ with index $b=779625$, say $E_1$ and $E_2$. By performing the commands ${\tt {N1:=PermutationGroup<299376|F(E1)>}}$ and ${\tt {N2:=PermutationGroup<299376|F(E2)>}}$, we obtain the subgroups of $N$ with index $b=779625$, namely $N_1$ and $N_2$. The 5 smallest orbit-lengths of $N_1$ and $N_2$ are 144, 288, $576^{3}$ and 48, $288^{2}$, $384^{2}$, respectively. But we have $k=480$ in case 51, which contradicts Lemma 2.3.

In case 61, we obtain $N$, the permutation representation of the group $J_{3}:2$ of degree 25840 and two subgroups $N_1$ and $N_2$ with index $b=147720$ by using the same methods as in case 51. The 5 smallest orbit-lengths of $N_1$ and $N_2$ are 64, $96^{3}$, 144 and 16, $48^{2}$, $72^{2}$, respectively. It is clear that there is no way to write $k=88$ as sum of some orbit-lengths of $N_1$, but for $N_2$, $88=16+72$. Suppose that $O_1$ is the orbit with length 16, and $O_4$ and  $O_5$ are the two orbits with length 72. Let $B_1=O_1\cup O_4$ and $B_2=O_1\cup O_5$. Then we have the block orbit-length $|B_1^{N}|=174420=|B_2^{N}|$. However, using the command ${\tt Design<2,25840|C>}$, where $C=B_1^{N}$ or $B_2^{N}$, we know that both structures are not $2$-designs.$\hfill\square$

\begin{lemma}\label{LEM3.7}{\rm}\quad
$G$ is not $M$.\end{lemma}
\textbf{Proof.}
Suppose $G$ is the Monster $M$. For $M$, the complete list of maximal subgroups has not been finished yet. At present, we obtain 43 known maximal subgroups of $M$ by inspecting the {\sc Atlas}, and we can check that none of them give rise to a set of parameters for $\mathcal{D}$ satisfying properties \rm (1)-\rm (5) in Lemma \ref{LEM2}. In $\cite{Bray}$, each maximal subgroup $H$ of $M$ which is not listed in the {\sc Atlas} is almost simple with Soc$(H)$ isomorphic to one of $L_2(13)$, $U_3(4)$, $U_3(8)$ and $Sz(8)$. They can also be ruled out easily by using the similar methods presented in subsection \ref{Sub3.1}. $\hfill\square$

\begin{lemma}\label{LEM3.8}{\rm}\quad
Theorem 1 holds in case 31.\end{lemma}
\textbf{Proof.}
In case 31, first we get the unique permutation representation of the Higman-Sims simple group $G=HS$ acting on 176 points by using ${\tt {G:=PrimitiveGroup(176,4)}}$. There are two conjugacy classes of subgroups with index 1100, denoted by $K_1$ and $K_2$ as representatives. Both $K_1$ and $K_2$ have two orbits in their action on the set of points. For $K_1$, the orbit-lengths are 8 and 168, and for $K_2$, the orbit-lengths are 56 and 120. There is no way to express $k=8$ as sum of some orbit-lengths of $K_2$. So we only need to consider $K_1$. Denote the orbits with lengths 8 and 168 by $\Delta_1$ and $\Delta_2$, respectively. It is easy to find that the block orbit-length $|\Delta_1^{G}|=1100$. Thus $\Delta_1$ can be regarded as a block and  the command ${\tt {Design<2,176|C>}}$ where $C$ refers to $\Delta_1^{G}$ returns a $2$-$(176,8,2)$ design with 1100 blocks.
 Let $\mathcal{P}=\{1, 2,\ldots, 176\}$. The orbits of $HS$ acting on $\mathcal{P}$ are listed as follows:
\begin{align*}
\Delta_1=\{ &18, 22, 24, 29, 57, 87, 166, 175 \},
\ \Delta_2=\mathcal{P}\backslash\Delta_1.
\end{align*}
The basic block of $\mathcal{D}$ is $\Delta_1$. Clearly, $HS \leq Aut(\mathcal{D})$ is block-transitive and point-primitive.$\hfill\square$

This completes the proof of Theorem \ref{thm1}.

\subsection*{Acknowledgements}
Thanks to Professor Sanming Zhou at University of Melbourne for corrections and some useful discussion which lead to the improvement of the paper.

\end{document}